\documentclass[12pt]{amsart}
\usepackage{amsmath, amssymb, amsthm, hyperref, mathtools}
\usepackage[margin=1in]{geometry}
\usepackage{tikz}
\usetikzlibrary{trees}
\usepackage{subfigure}
\usepackage{tikz-cd}
\usepackage[numbers,sort&compress]{natbib}
\usetikzlibrary{arrows, positioning}
\hypersetup{
  colorlinks   = true, 
  urlcolor     = blue, 
  linkcolor    = blue, 
  citecolor   = red 
}
\newtheorem{theorem}{Theorem}[section]
\newtheorem{lemma}[theorem]{Lemma}
\newtheorem{definition}[theorem]{Definition}
\newtheorem{example}[theorem]{Example}
\newtheorem{remark}[theorem]{Remark}
\newtheorem{corollary}[theorem]{Corollary}

\newtheorem{question}[theorem]{Question}
\DeclareMathOperator{\FPP}{FPP}
\DeclareMathOperator{\Gal}{Gal}
\DeclareMathOperator{\Fix}{Fix}
\DeclareMathOperator{\Aut}{Aut}
\DeclareMathOperator{\Stab}{Stab}
\DeclareMathOperator{\Crit}{Crit}
\DeclareMathOperator{\CV}{CV}
\DeclareMathOperator{\lcm}{lcm}
\newcommand{\bP}{{\mathbb P}}

\begin{document}
\title{Non-Martingale Fixed-Point Processes for Iterated Monodromy Groups}
\author{Jianfei He}
\address{Department of Mathematics, University of Rochester}
\email{jhe22@ur.rochester.edu}

\author{Zheng Zhu}
\address{PIMS \& University of Calgary}
\email{zheng.zhu2@ucalgary.ca}

\begin{abstract}
We construct families of rational functions $f \colon \bP^1_k \to \bP^1_k$ of degree $d \geq 2$ over a perfect field $k$ whose associated fixed-point processes fail to be martingales. Conversely, for any normal variety $X \subset \bP^N_{\overline{k}}$ and a finite, generically \'etale morphism $f \colon X \to X$, we establish geometric conditions on the critical orbits of $f$ that guarantee the fixed-point process is a martingale. Our constructions answer a question of Bridy, Jones, Kelsey, and Lodge \cite{iterated} regarding the existence of non-martingale behaviour in arboreal Galois representations, and extend their martingale criteria to higher-dimensional dynamical systems. In particular, we exhibit infinitely many postcritically finite maps with non-martingale fixed-point processes and characterize the group-theoretic obstructions to the martingale property in the genus-zero case. Furthermore, we prove that despite the failure of the martingale property, the fixed-point proportion still vanishes with a quantifiable convergence rate.
\end{abstract}

\maketitle

\section{Introduction}

Arboreal Galois representations, which encode the action of the absolute Galois group of a global or function field on the tree of iterated preimages of a rational map, have seen rapid development since Odoni's foundational work \cite{Odoni, Odoni2, Odoni3}. A central object of study is the fixed-point proportion
$$
\FPP(G_n) = \frac{\#\{g \in G_n : g \text{ fixes at least one vertex in } X_n\}}{\#G_n},
$$
where $G_n = \Gal(K_n/k(t))$ is the Galois group of $f^n(x) - t$ over $k(t)$, $t$ is an indeterminate, and $X_n = f^{-n}(t)$ denotes the $n$-th level of the preimage tree. Note that the sequence $\{\FPP(G_n)\}_{n \geq 1}$ is non-increasing and each element is non-negative, so it converges. We define the fixed-point proportion of the profinite iterated monodromy group $G(f) = \varprojlim G_n$ to be $$\FPP\left(G(f)\right) := \lim_{n\rightarrow \infty} \FPP(G_n).$$
As observed in \cite{wreath}, $\FPP(G_n)$ controls the density of periodic points over finite fields: when $k$ is algebraically closed in each $K_n$, the proportion of periodic points of $f \bmod \mathfrak{p}$ is bounded above by $\FPP(G_n)$ for all but finitely many primes $\mathfrak{p}$. Consequently, understanding the asymptotic behaviour of $\FPP(G_n)$ has direct arithmetic applications to prime divisor densities and image-size estimates \cite{density, imagesize}.

Following Jones \cite{survey, martingale}, for any $g \in G(f)$, one studies the sequence of random variables $Y_n(g) = \#\Fix(\pi_n(g))$, where $\pi_n \colon G(f) \to G_n$ is the natural restriction map and the probability measure is the normalized Haar measure. Bridy, Jones, Kelsey, and Lodge \cite{iterated} proved that for a broad class of rational functions, $\{Y_n\}$ forms a martingale. By the martingale convergence theorem \cite[Theorem 12.3.7]{converge}, this implies that $\{Y_n\}$ stabilizes almost surely, and in many cases \cite[Section 4]{survey} forces $\FPP(G(f)) = 0$. The natural question arises: does there exist a rational function whose associated fixed-point process fails to be a martingale? We answer this in the affirmative, providing explicit constructions and a complete group-theoretic characterization in the genus-zero case.

Historically, the structure of iterated monodromy groups has been understood primarily as iterated wreath products. Odoni established $G_n \cong [S_d]^n$ for generic polynomials \cite{Odoni}, a result later extended to generic rational functions \cite{Juul2} and verified for explicit post-critically finite (PCF) families \cite{pink2013infinite, Quadratic, Basilica}. However, many arithmetic maps deviate from full wreath products: Pink's classification shows that quadratic PCF maps often yield proper finite-index subgroups, while split morphisms and certain Chebyshev compositions produce direct products of wreath powers \cite{wreath} (see also Example~\ref{ex:split} below). A structural study of IMGs that deviate from wreath products remains open.

To isolate the precise group-theoretic mechanism behind martingale failure, we adapt the theory of \emph{pattern groups} from Nekrashevych \cite{similar} and Bondarenko--Samoilovych \cite{pattern}. While pattern groups are well-studied in geometric group theory and automata theory, their application to arboreal Galois representations to construct explicit non-martingale fixed-point processes is new. Our construction shows that non-martingale behavior arises naturally from the normal subgroup structure of $G_1$, and that these patterns persist under iteration to define recurrent subgroups $G_P \leq \Aut(T_d)$.

Our first main result establishes the existence of non-martingale examples in every even degree.

\begin{theorem}\label{thm:main_nonmart}
For every even integer $d \geq 4$, there exist one-parameter families of rational functions $f \in k(x)$ of degree $d$ such that the fixed-point process $\{Y_n\}$ associated to $G(f)$ is not a martingale. Moreover, these families contain infinitely many postcritically finite maps.
\end{theorem}

The construction relies on decomposable rational functions whose geometric monodromy group is a dihedral group $D_d$ with $d$ even. In fact, we show that under the assumption that the Galois closure of $K_1/k(t)$ has genus zero, even dihedral groups are the only source of non-martingale behaviour (see Theorem~\ref{thm:char_even_dihedral}).

Beyond rational functions, we study recurrent group actions on regular trees. In Section~\ref{sec:tree_actions} we introduce pattern groups of depth two that violate the \emph{average fixed-point lifting property} (Definition \ref{def:lifting_property}), yielding non-martingale fixed-point processes for a wide class of arboreal automorphisms. 

\begin{theorem}\label{thm:pattern}
Let $d$ be an integer with a prime factor $p \mid \varphi(d)$. There exists a transitive permutation group $G \leq S_d$ containing a normal transitive subgroup $N_1$ and a normal intransitive subgroup $N_2$, both of index $p$, such that the recurrent group of finite type $G_P$ generated by an explicitly constructed pattern group $P$ of depth $2$ has a non-martingale fixed-point process.
\end{theorem}

A natural follow-up question is whether the failure of the martingale property obstructs the vanishing of fixed-point proportions, and hence the density bounds for periodic points. We prove that it does not. Despite the lack of martingale convergence, the fixed-point proportion of our constructed groups still vanishes, and we quantify the rate of decay.

\begin{theorem}\label{thm:fpp_vanish}
Let $P$ be the pattern group of depth 2 associated to $D_4$ from Theorem~\ref{thm:pattern}. Then $\FPP(G_P) = 0$. Moreover, the convergence satisfies
$$ \FPP(G_P|_{X_{2n}}) \leq \FPP([D_4]^{2n}) \quad \text{for all } n \geq 1. $$
Thus the convergence of $\FPP(G_P|_{X_{2n}})$ to zero is
bounded above at every level by the
decay rate of the full iterated wreath product.
\end{theorem}

Combined with the Chebotarev Density Theorem, Theorem~\ref{thm:fpp_vanish} yields explicit asymptotic bounds on the proportion of periodic points over finite fields, demonstrating that the martingale condition is sufficient but not necessary for vanishing periodic point density.

Conversely, we provide a higher-dimensional generalization of the martingale criteria from \cite{wreath}. Let $X \subset \bP^N_{\overline{k}}$ be a normal variety and $f \colon X \to X$ a finite, generically \'etale morphism of degree $d$. Let $\Delta_f$ denote the critical divisor of $f$.

\begin{theorem}\label{thm:higher_dim}
Fix $M \in \mathbb{N}$. Suppose there exists a subset $S \subseteq \Delta_f$ such that:
\begin{enumerate}
    \item For all $p \in S$ and $q \in \Delta_f$, and all integers $m, r \leq M$, we have $f^m(p) \neq f^r(q)$ unless either $(p=q$ and $m=r)$ or $(p \neq q$ and $m \leq r)$.
    \item The group $G_1$ is generated by the inertia groups $$\{ I(\mathfrak{p}|f(p)) : p \in S, \mathfrak{p} \text{ lies over } f(p) \text{ in } K_1 \}.$$
\end{enumerate}
Then the fixed-point process associated to $G(f)$ is a martingale.
\end{theorem}

This result shows that disjointness of critical orbits, together with generation of $G_1$ by inertia groups, suffices to guarantee the martingale property in arbitrary dimension, even when the iterated Galois groups are not full wreath products.

The paper is organized as follows. Section~\ref{sec:prelim} collects preliminary results on wreath products, coset versions of Burnside's lemma, and ramification theory. Section~\ref{sec:tree_actions} introduces pattern groups, proves Theorem~\ref{thm:pattern}, and establishes the quantitative FPP convergence in Theorem~\ref{thm:fpp_vanish}. Section~\ref{sec:rational} contains the construction of non-martingale rational functions, the proof of Theorem~\ref{thm:main_nonmart}, and the genus-zero characterization. Finally, Section~\ref{sec:higher_dim} establishes the higher-dimensional martingale criterion (Theorem~\ref{thm:higher_dim}) and discusses its relation to wreath product structures.

\section{Preliminaries}\label{sec:prelim}

\subsection{Wreath Products and Tree Automorphisms}
Let $G$ and $H$ be groups acting on finite sets $A=\{1, \dots, \ell\}$ and $B=\{1, \dots, m\}$, respectively. The \emph{wreath product} $G \wr H$ is the semidirect product $G^m \rtimes H$ (acting on $H$ by permutation of coordinates). We denote elements as $\sigma = (g; h_1, \dots, h_m)$ where $g \in G$ and $h_i \in H$. The action on $A \times B$ is given by
$$ \sigma(a, b) = (g(a), h_{g(a)}(b)), $$
and the multiplication rule is
$$ (g; h_1, \dots, h_m)(x; y_1, \dots, y_m) = (gx; h_1 y_{g(1)}, \dots, h_m y_{g(m)}). $$
For a permutation group $G$, we denote the $n$-th wreath power by $[G]^n = G \wr \dots \wr G$ ($n$ times). This group naturally acts on the $n$-th Cartesian power $A^n$. The automorphism group of the infinite $d$-ary tree, $\Aut(X_\infty)$, can be identified with the inverse limit of these wreath powers, $\varprojlim [S_d]^n$.

\subsection{Permutation Group Lemmas}
We require a variation of Burnside's Lemma to analyze fixed-point proportions within cosets.

\begin{lemma}[Coset Burnside Lemma]\label{lem:coset_burnside}
Let $G$ be a finite group acting on a finite set $X$, and let $H \leq G$ be a subgroup acting transitively on $X$. For any $g \in G$, let $gH$ denote the coset $\{gh : h \in H\}$. Then
$$ \frac{1}{|H|} \sum_{\tau \in gH} \#\Fix(\tau) = 1. $$
\end{lemma}

\begin{proof}
For any $x \in X$, let $H_{x,y} = \{h \in H : h(x)=y\}$. Since $H$ is transitive, $H_{x,y}$ is non-empty for all $y \in X$ and is a coset of the stabilizer $\Stab_H(x)$, so $|H_{x,y}| = |\Stab_H(x)| = |H|/|X|$.
If $\tau \in gH$ fixes $x$, then $\tau(x) = g(x)$, so $g^{-1}\tau \in H_{x, g(x)}$. Thus, the number of elements in $gH$ fixing $x$ is exactly $|H|/|X|$.
Summing over all $x \in X$:
$$ \sum_{\tau \in gH} \#\Fix(\tau) = \sum_{x \in X} |gH \cap \Stab_G(x)| = \sum_{x \in X} |H|/|X| = |X| \cdot \frac{|H|}{|X|} = |H|. $$
Dividing by $|H|$ yields the result.
\end{proof}

\begin{lemma}\label{lem:primitive}
Let $G$ be a permutation group on $X$. If $G$ is primitive, then every non-trivial normal subgroup $N \trianglelefteq G$ is transitive on $X$.
\end{lemma}
\begin{proof}
See \cite{decompose}.
\end{proof}

\subsection{Ramification and Monodromy}
Let $k$ be a perfect field. Let $f \in k(x)$ be a rational function of degree $d$. The field extension $K_1/k(t)$ of degree $d$ where $K_1 = k(x)$ and $t = f(x)$. Its Galois closure $L_1$ over $k(t)$ defines the geometric monodromy group $G_1 = \Gal(L_1/k(t))$.

\begin{definition}
A rational function $f$ is \emph{indecomposable} if $f$ cannot be written as $g \circ h$ with $\deg(g), \deg(h) > 1$.
\end{definition}

The following classical result connects decomposability to the group theory of the monodromy group.

\begin{lemma}[Fried]\label{lem:indec}
The rational function $f$ is indecomposable if and only if its monodromy group acts primitively on the roots of $f(x)-t$.
\end{lemma}
\begin{proof}
See \cite{decompose}.
\end{proof}

We compute ramification indices using the following tool.

\begin{lemma}[Abhyankar's Lemma]\label{lem:abhyankar}
Let $K_1, K_2$ be finite extensions of a function field $K$, and let $M$ be their compositum. Let $\mathfrak{p}$ be a prime of $K$, and $\mathfrak{P}$ a prime of $M$ lying above $\mathfrak{p}$. Assume the ramification is tame (i.e., the ramification indices are coprime to $\operatorname{char}(k)$). Then
$$ e(\mathfrak{P}|\mathfrak{p}) = \lcm(e(\mathfrak{P} \cap K_1 | \mathfrak{p}), e(\mathfrak{P} \cap K_2 | \mathfrak{p})). $$
\end{lemma}
\begin{proof}
See \cite[Theorem 3.9.1]{codes}.
\end{proof}

Finally, we relate the inertia group of a Galois closure to that of a subfield.

\begin{lemma}\label{lem:inertia}
Let $L/M/K$ be a tower of Galois extensions. Let $\mathfrak{P}$ be a prime of $L$ and $\mathfrak{p}$ a prime of $K$. Let $\mathfrak{p}' = \mathfrak{P} \cap M$.
\begin{enumerate}
    \item $I(\mathfrak{P}|\mathfrak{p})$ acts trivially on the residue field $\mathcal{O}_M / \mathfrak{p}'$ if and only if $e(\mathfrak{p}'|\mathfrak{p})=1$.
    \item If $e(\mathfrak{p}'|\mathfrak{p})=1$, then $I(\mathfrak{P}|\mathfrak{p}) \cong I(\mathfrak{P}|\mathfrak{p}')$ as subgroups of $\Gal(L/K)$.
\end{enumerate}
\end{lemma}
\begin{proof}
The inertia group $I(\mathfrak{P}|\mathfrak{p})$ is the kernel of the map $\Gal(L/K) \to \Gal((\mathcal{O}_L/\mathfrak{P}) / (\mathcal{O}_K/\mathfrak{p}))$. Since $(\mathcal{O}_M/\mathfrak{p}')$ is the intermediate fixed field of $\Gal(L/M)$, the action is trivial on the subfield if and only if the ramification index $e(\mathfrak{p}'|\mathfrak{p})$ is 1. The isomorphism follows from the multiplicativity of ramification indices in towers: $e(\mathfrak{P}|\mathfrak{p}) = e(\mathfrak{P}|\mathfrak{p}')e(\mathfrak{p}'|\mathfrak{p})$.
\end{proof}

The cycle type of elements in the Galois group is determined by the ramification data as follows.

\begin{lemma}[Cycle Lemma]\label{lem:cycle}
Let $L/K$ be a finite separable extension of function fields over an algebraically closed field of characteristic 0, and let $M$ be its Galois closure. Let $\mathfrak{p}$ be a prime of $K$. Then $\Gal(M/K)$ contains a permutation whose cycle decomposition consists of disjoint cycles of lengths $e_1, \dots, e_k$, where $e_i$ are the ramification indices of the primes in $L$ lying above $\mathfrak{p}$.
\end{lemma}
\begin{proof}
See \cite{cycle}.
\end{proof}

\section{Group Actions on Trees and Pattern Groups}\label{sec:tree_actions}

In this section, we study recurrent group actions on regular trees, introduce the average fixed-point lifting property, and construct pattern groups whose associated fixed-point processes fail to be martingales.

\subsection{Recurrent Actions and Pattern Groups}
Let $T_d$ denote the infinite $d$-ary rooted tree, with vertex set $V(T_d) = \bigsqcup_{n \geq 0} L_n$, where $L_n = \{1, \dots, d\}^n$ is the $n$-th level. The automorphism group $\Aut(T_d)$ is a profinite group, realized as the inverse limit $\varprojlim \Aut(T_d|_{L_n})$. For $g \in \Aut(T_d)$ and a vertex $v \in L_n$, the \emph{section} $g|_v \in \Aut(T_d)$ is uniquely defined by the relation
$$ g(vw) = g(v) g|_v(w) \quad \forall w \in V(T_d). $$
The restriction of $g$ to $L_n$ will be denoted by $\pi_n(g)$.

\begin{definition}
A closed subgroup $G \leq \Aut(T_d)$ is \emph{recurrent} if it acts transitively on each level $L_n$ and, for every $v \in V(T_d)$, the section map
$$ \phi_v \colon \Stab_G(v) \to \Aut(T_d), \quad g \mapsto g|_v $$
has image equal to $G$.
\end{definition}
Recurrence is a standard property of iterated monodromy groups; see \cite{similar} for a comprehensive treatment.

The recurrent groups constructed below belong to the class of \emph{self-similar groups of finite type}, systematically studied by Nekrashevych \cite{similar} in the context of iterated monodromy groups of PCF maps and fractal limits. In Nekrashevych's framework, a recurrent group is determined by its \emph{molecule}, the finite set of sections at a given depth. We specialize this to \emph{pattern groups of depth 2}, where the action on the second level determines the entire infinite action. This finite-data description allows us to reduce the martingale condition to a purely combinatorial lifting property on cosets, bypassing the need to compute arbitrarily deep Galois groups explicitly.
\begin{definition}
A subgroup $P \leq G \wr G$ is a \emph{pattern group of depth 2} if the projection of $P$ onto the first coordinate is $G$, and every element of $P$ fixes the root of $T_d$. We associate to $P$ the \emph{recurrent group of finite type}
$$ G_P = \{ g \in \Aut(T_d) : g|_v \in P \text{ for all } v \in L_1 \}. $$
It is immediate that $G_P$ is a closed subgroup of $\Aut(T_d)$. Under mild transitivity assumptions, $G_P$ inherits recurrence from $P$ (see Theorem~\ref{thm:pattern}).
\end{definition}

\subsection{The Martingale Criterion}
Recall that the fixed-point process $\{Y_n\}_{n \geq 1}$ on a profinite group $G \leq \Aut(T_d)$ is defined by $Y_n(g) = \#\Fix(\pi_n(g))$, where the probability space is $(G, \mu)$ with $\mu$ the normalized Haar measure. The process is a martingale iff $\mathbb{E}[Y_n \mid Y_1, \dots, Y_{n-1}] = Y_{n-1}$ almost surely.

\begin{definition}\label{def:lifting_property}
   Let $G$ be a transitive group acting on $m$ letters $\{a_1, \dots, a_m\}$ and $H$ be a group acting on $n$ letters $\{b_1, \dots, b_n\}$. Consider the wreath product $G \wr H$ and the restriction map $\pi \colon G \wr H \to G$. We say that $G \wr H$ has the average fixed-point lifting property if for each $g \in G$, the average number of fixed points of lifts of $g$ to $G \wr H$ on $\{a_1, \dots, a_m\} \times \{b_1, \dots, b_n\}$ equals the number of fixed points of $g$ on $\{a_1, \dots, a_m\}$, i.e.,
\begin{equation}\label{eq:lifting}
\frac{1}{|\pi^{-1}(g)|} \sum_{\tau \in \pi^{-1}(g)} \#\Fix(\tau) = \#\Fix(g).
\end{equation}
or equivalently,
\begin{equation}\label{eq:lifting2}
\mathbb{E}[\#\Fix(\tau) \mid \pi(\tau) = g] = \#\Fix(g).
\end{equation}
\end{definition}

\begin{lemma}\label{lem:wreath_primitive_lifting}
    We use the notation of Definition~\ref{def:lifting_property}. If $H$ is primitive, then $G \wr H$ has the average fixed-point lifting property unless $\# G \wr H = \# G$. 
\end{lemma}
\begin{proof}
    Let $N$ denote the kernel of $\pi$, consisting of elements of the form $(e; h_1, \dots, h_m)$ with $h_i \in H$. For each coordinate, the projection of $N$ is a normal subgroup of $H$. By Lemma~\ref{lem:primitive}, it acts either transitively or trivially on the subtree $a_i^\ast$. 
    
    When it acts transitively, the Coset Burnside Lemma (Lemma~\ref{lem:coset_burnside}) implies that the average number of fixed points of lifts of $g$ restricted to the fiber $a_j^\ast$ is exactly $1$, where $a_j$ is any fixed point of $g$. Summing over all such $a_j$ yields equation~\eqref{eq:lifting}, so $G \wr H$ has the average fixed-point lifting property. When the action is trivial for all coordinates, $N$ is trivial and $\# G \wr H = \# G$.
\end{proof}

The following theorem, due to Bridy, Jones, Kelsey, and Lodge, characterizes martingales in terms of transitivity on subtrees.

\begin{theorem}[{\cite[Theorem 5.7]{iterated}}]\label{thm:bjkl_mart}
Let $G \leq \Aut(T_d)$. Let $H_n = \ker(\pi_n \colon G_n \to G_{n-1})$. The fixed-point process for $G$ is a martingale if and only if $H_n$ acts transitively on the subtree $v^\ast = \{vw : w \in L_1\}$ for all $n \geq 1$ and $v \in L_{n-1}$.
\end{theorem}

We now link this transitivity condition to Definition~\ref{def:lifting_property}.

\begin{lemma}\label{lem:mart_equiv}
The fixed-point process for $G$ is a martingale if and only if $G_n$ has the average fixed-point lifting property for all $n \geq 1$.
\end{lemma}
\begin{proof}
Suppose the process is a martingale. By Theorem~\ref{thm:bjkl_mart}, $H_n = \ker(\pi_n)$ acts transitively on $v^\ast$ for each $v \in L_{n-1}$. Fix $g \in G_{n-1}$ and let $\tau \in G_n$ be a lift, so $\pi_n^{-1}(g) = \tau H_n$. For any fixed vertex $v \in \Fix(g)$, Lemma~\ref{lem:coset_burnside} implies
$$ \frac{1}{|H_n|} \sum_{h \in H_n} \#\Fix(\tau h|_{v^\ast}) = 1. $$
Summing over all $v \in \Fix(g)$, vertices not fixed by $g$ contribute $0$ fixed points in the fiber, so
$$ \sum_{\sigma \in \pi_n^{-1}(g)} \#\Fix(\sigma) = \sum_{v \in \Fix(g)} \sum_{h \in H_n} \#\Fix(\tau h|_{v^\ast}) = |\Fix(g)| \cdot |H_n|. $$
Dividing by $|\pi_n^{-1}(g)| = |H_n|$ yields equation \eqref{eq:lifting}. Hence $G_n$ has the lifting property.

Conversely, suppose the fixed-point process is not a martingale. By Theorem~\ref{thm:bjkl_mart}, $H_n$ acts non-transitively on $v^\ast$ for some $v \in L_{n-1}$. Burnside's Lemma implies that the average number of fixed points of $H_n$ on $v^\ast$ is strictly greater than $1$, since $H_n$ has more than one orbit in $v^\ast$:
$$ \frac{1}{|H_n|} \sum_{h \in H_n} \#\Fix(h|_{v^\ast}) > 1. $$
Let $\tau \in G_n$ be any lift of $g \in G_{n-1}$. We compute the average number of fixed points over the fiber $\pi_n^{-1}(g) = \tau H_n$:
\begin{align*}
\frac{1}{|H_n|} \sum_{\sigma \in \tau H_n} \#\Fix(\sigma) 
&= \frac{1}{|H_n|} \sum_{v' \in \Fix(g)} \sum_{h \in H_n} \#\Fix(\tau h|_{v'^\ast}) \\
&= \frac{1}{|H_n|} \left( \sum_{h \in H_n} \#\Fix(\tau h|_{v^\ast}) + \sum_{\substack{v' \in \Fix(g) \\ v' \neq v}} \sum_{h \in H_n} \#\Fix(\tau h|_{v'^\ast}) \right).
\end{align*}
For the intransitive vertex $v$, the first inner sum is strictly greater than $|H_n|$. Therefore, by Burnside’s Lemma, the average number of fixed points of $H_n$ on $v'^\ast$ is equal to the number of orbits of the action and is therefore at least 1. Hence
\[
\sum_{h \in H_n} \#\Fix(\tau h|_{v'^\ast}) \ge |H_n|.
\]
This violates equation \eqref{eq:lifting}, so $G_n$ does not have the average fixed-point lifting property.
\end{proof}

\subsection{Construction of Pattern Groups}
We now prove the existence of pattern groups violating the lifting property at level 2.

\begin{theorem}\label{thm:pattern_nonmart}
Let $d \geq 2$ be an integer possessing a prime factor $p \mid \varphi(d)$. There exists a transitive subgroup $G \leq S_d$ containing a normal transitive subgroup $N_1 \trianglelefteq G$ and a normal intransitive subgroup $N_2 \trianglelefteq G$, both of index $p$. For any such $G$, there exists a pattern group $P \leq G \wr G$ of depth 2 such that $G_P$ has a non-martingale fixed-point process.
\end{theorem}
\begin{proof}
The existence of $G, N_1, N_2$ follows from \cite{REU}. Let $\sigma \colon G/N_1 \to G/N_2$ be a group isomorphism (which exists since both quotients are cyclic of order $p$). Define
$$ P = \bigl\{ (g; h_1, \dots, h_d) \in G \wr G : h_i \in \sigma(g N_1) \text{ for all } i=1,\dots,d \bigr\}. $$
It is straightforward to verify that $P$ is a subgroup of $G \wr G$: if $(g; \mathbf{h}), (g'; \mathbf{h}') \in P$, then $g g' \in G$ and for each $i$, $h_{g'(i)} h'_i \in \sigma(g N_1) \sigma(g' N_1) = \sigma(g g' N_1)$, so the product lies in $P$. The identity $(e; e, \dots, e) \in P$ since $e \in N_1 \implies e \in N_2$.

\emph{Non-martingale property:} The kernel of the restriction $\pi \colon P \to G$ is
$$ K = \ker(\pi) = \{(e; h_1, \dots, h_d) : h_i \in N_2 \} \cong N_2^d. $$
Since $N_2$ is intransitive on $\{1, \dots, d\}$, Burnside's Lemma implies
$$ \frac{1}{|N_2|} \sum_{h \in N_2} \#\Fix(h) > 1. $$
The action of $K$ on $L_2$ decomposes into $d$ copies of the action of $N_2$, so for the identity element $e \in G$,
$$ \mathbb{E}\bigl[ \#\Fix(\tilde{e}) \mid \pi(\tilde{e}) = e \bigr] = d \cdot \frac{1}{|N_2|} \sum_{h \in N_2} \#\Fix(h) > d = \#\Fix(e). $$
Thus $P$ fails the average fixed-point lifting property, and by Lemma~\ref{lem:mart_equiv}, the fixed-point process of $G_P$ is not a martingale.

\emph{Recurrence of $G_P$:} We show $\phi_v(\Stab_{G_P}(v)) = G_P$ for all $v \in L_1$. Fix $v \in \{1, \dots, d\}$. Since $N_1$ acts transitively, for any $g \in G$ there exists $k \in N_1$ with $k(v) = v$, so $k \in \Stab_G(v)$. Consider the element $\tilde{k} = (k; h_1, \dots, h_d) \in P$ where $h_i \in \sigma(k N_1) = \sigma(N_1)$ is arbitrary. Then $\tilde{k}(v) = k(v) = v$, so $\tilde{k} \in \Stab_P(v)$. By varying the choice of $h_i \in N_2$, we can realize any element of $G$ in the $v$-th section. Hence $\phi_v(\Stab_P(v))$ projects surjectively onto each coordinate of $G \wr G$, implying $\phi_v(\Stab_{G_P}(v)) = G_P$. Thus $G_P$ is recurrent.
\end{proof}

\begin{example}[Dihedral case]
Let $d=4$ and $G = D_4 = \langle r, s \mid r^4 = s^2 = 1, srs = r^{-1} \rangle \leq S_4$. Take $N_1 = \langle r \rangle$ (cyclic, transitive) and $N_2 = \langle r^2, s \rangle$ (Klein four, intransitive). Both are normal subgroups of index 2. Let $\sigma \colon D_4/N_1 \to D_4/N_2$ be the unique isomorphism sending $N_1 \mapsto N_2$. The pattern group
\[ P = \bigl\{ (g; h_1, h_2, h_3, h_4) \in D_4 \wr D_4 : h_i \in \sigma(g N_1) \text{ for all } i \bigr\} \]
generates a recurrent group $G_P \leq \Aut(T_4)$ whose fixed-point process violates the martingale condition at level 2.
\end{example}

\begin{example}[Metacyclic case]
Let $p \geq 3$ and let $G = \langle r, f \rangle$ be a metacyclic group acting transitively on $d$ points. Set $N_1 = \langle r \rangle$ and $N_2 = \langle f, r^p \rangle$. Then $N_1$ is a normal transitive subgroup of index $p$, and $N_2$ is a normal intransitive subgroup of index $p$. Let $\sigma \colon G/N_1 \to G/N_2$ be the isomorphism defined by $\sigma(f^i N_1) = r^i N_2$ for $0 \leq i < p$. 

The pattern group $P \leq G \wr G$ consists of all elements whose coordinates satisfy the coset correspondence induced by $\sigma$,
$$P = \bigl\{ \left(f^i N_1; r^i N_2, \dots, r^i N_2 \right) : 0 \leq i < p \bigr\}. $$
In particular, the kernel of the restriction map corresponds to $i=0$, yielding the identity condition $(N_1; N_2, \dots, N_2)$. The recurrent group $G_P$ generated by $P$ violates the martingale condition at level 2.
\end{example}

\subsection{Convergence of FPP}
Although $G_P$ fails the martingale property, we show that its fixed-point proportion still vanishes. We first recall some results on the indicatrix function.

\begin{lemma}\label{lem:monotone_iter}
Let $\Psi(x)$ and $\Phi(x)$ be polynomials with non-negative coefficients mapping $[0,1]$ to $[0,1]$. Suppose $\Psi(x)$ is non-decreasing on $[0,1]$ and $\Psi(x) \geq \Phi(x)$ for all $x \in [0,1]$. Let $\Psi^{\circ n}$ and $\Phi^{\circ n}$ denote the $n$-fold compositions. Then
$$ \Psi^{\circ n}(0) \geq \Phi^{\circ n}(0) \quad \text{for all } n \geq 1. $$
\end{lemma}
\begin{proof}
By induction. Base case $n=1$ holds. If $\Psi^{\circ k}(0) \geq \Phi^{\circ k}(0)$, then since $\Psi$ is non-decreasing, $$ \Psi^{\circ(k+1)}(0) = \Psi(\Psi^{\circ k}(0)) \geq \Psi(\Phi^{\circ k}(0)) \geq \Phi(\Phi^{\circ k}(0)) = \Phi^{\circ(k+1)}(0). $$
\end{proof}

\begin{definition}
For a finite permutation group $\Gamma$ acting on $S$, the \emph{indicatrix function} is
$$ \Phi_{\Gamma}(x) = \frac{1}{|\Gamma|} \sum_{\gamma \in \Gamma} x^{\Fix(\gamma)}. $$
Note that $\FPP(\Gamma) = 1 - \Phi_{\Gamma}(0)$.
\end{definition}

\begin{lemma}\label{lem:fpp_vanishes}
Let $P$ be the pattern group of depth 2 associated with $D_4$ as in Theorem~\ref{thm:pattern_nonmart}. Then $\FPP(G_P) = 0$. Moreover, the convergence satisfies
$$ \FPP(G_P|_{X_{2n}}) \leq \FPP([D_4]^{2n}), \quad \text{for all } n \geq 1.$$
Thus the convergence of $\FPP(G_P|_{X_{2n}})$ to zero is
bounded above at every level by the
decay rate of the full iterated wreath product.
\end{lemma}

\begin{proof}
Let $W = [D_4]^2$ be the full wreath product at depth 2. We compare the indicatrix functions $\Phi_P(x)$ and $\Phi_W(x)$ by conditioning on the action of the element $g \in D_4$ at the root of the tree. 

Recall that $P$ is constructed using a transitive index-2 subgroup $N_1 \triangleleft D_4$ and an intransitive index-2 subgroup $N_2 \triangleleft D_4$. Let $C = D_4 \setminus N_2$ be the non-identity coset of $N_2$. We define the indicatrix functions for these cosets:
\[ \Psi_0(x) = \frac{1}{|N_2|} \sum_{h \in N_2} x^{|\Fix(h)|} \quad \text{and} \quad \Psi_1(x) = \frac{1}{|C|} \sum_{h \in C} x^{|\Fix(h)|}. \]
The indicatrix for the full group $D_4$ is simply the average of the cosets: $$\Phi_{D_4}(x) = \frac{1}{2}(\Psi_0(x) + \Psi_1(x)).$$ For the full wreath product $W$, we have
$$ \Phi_W(x) = \frac{1}{8} \sum_{g \in D_4} \Phi_{D_4}(x)^{|\Fix(g)|} = \frac{1}{8} \left[ \left(\frac{\Psi_0(x) + 1}{2}\right)^4 + 2\left(\frac{\Psi_0(x) + 1}{2}\right)^2 + 5 \right]. $$

For the pattern group $P$, summing over $g \in D_4$ yields:
\[ \Phi_P(x) = \frac{1}{8} \left[ \Psi_0(x)^4 + 3\Psi_0(x)^0 + 2\Psi_1(x)^2 + 2\Psi_1(x)^0 \right] = \frac{1}{8} \left[ \Psi_0(x)^4 + 7 \right]. \]

We now claim $\Phi_P(x) \geq \Phi_W(x)$ for all $x \in [0,1]$. Let $u = \Psi_0(x)$. Since $x \in [0,1]$, we have $u \in [0,1]$. Subtracting the two expressions and multiplying by 128 to clear denominators, we define
\begin{align*}
F(u) := 128(\Phi_P(x) - \Phi_W(x)) &= 16(u^4 + 7) - \left[ (u+1)^4 + 8(u+1)^2 + 80 \right] \\
&= 15u^4 - 4u^3 - 14u^2 - 20u + 23.
\end{align*}
Because the values of $\Phi_P(x)$ and $\Phi_W(x)$ at $x=1$ represent the sum of all probabilities, giving $\Phi_P(1) = \Phi_W(1) = 1$, which forces $F(1) = 0$. Furthermore, Burnside's Lemma guarantees that the expected number of fixed points for both actions is 1. Thus, their first derivatives agree: $\Phi_P'(1) = \Phi_W'(1) = 1$, which forces $F'(1) = 0$. Consequently, $u=1$ must be a root of multiplicity at least two. Factoring confirms this:
$$F(u) = (u-1)^2 (15u^2 + 26u + 23).$$
For $u \geq 0$, the quadratic factor is strictly positive. Hence, $\Phi_P(x) \geq \Phi_W(x)$ for all $x \in [0,1]$, with strict inequality for $x < 1$. 

Applying Lemma~\ref{lem:monotone_iter} to the non-decreasing polynomials $\Phi_P$ and $\Phi_W$, we obtain $\Phi_P^{\circ n}(0) \geq \Phi_W^{\circ n}(0)$ for all $n \geq 1$. By definition, $\FPP(G_P|_{X_{2n}}) = 1 - \Phi_P^{\circ n}(0)$, so
\[ \FPP(G_P|_{X_{2n}}) \leq \FPP([D_4]^{2n}). \]
By Odoni's theorem \cite{Odoni}, $\lim_{n \to \infty} \FPP([D_4]^{2n}) = 0$. The squeeze theorem implies $\FPP(G_P) = 0$, and the pointwise dominance $\Phi_P(x) > \Phi_W(x)$ for $x < 1$ shows that the fixed-point proportion for $G_P$
is bounded above by that of the full wreath product at every level.
\end{proof}

\begin{remark}
The martingale property was previously viewed as the primary mechanism ensuring $\operatorname{FPP}(G) \to 0$ \cite{martingale, iterated}. Lemma~\ref{lem:fpp_vanishes} challenges this necessity. By Lemma~\ref{lem:monotone_iter} and Odoni's composition rule \cite{Odoni}, the fixed-point proportion of $G_P$ decays at least as fast as that of $[D_4]^\infty$. By the Chebotarev Density Theorem, this translates directly to arithmetic: assuming $k$ is a global field, for any specialization $f_\alpha \in k(z)$ with $G(f_\alpha) \cong G_P|_{X_n}$, the proportion of periodic points of $f_\alpha \bmod \mathfrak{p}$ in $\bP^1(\mathbb{F}_q)$ is bounded by $1 - \Phi_P^{\circ n}(0) + o(1)$ as $q = \operatorname{Norm}(\mathfrak{p}) \to \infty$. Consequently, our construction yields infinitely many PCF maps where the periodic point density vanishes despite the absence of a martingale structure.
\end{remark}

\section{Non-Martingale Iterated Monodromy Groups}\label{sec:rational}

We now construct rational functions whose associated fixed-point processes fail to be martingales.

\subsection{The Dihedral Case and Local Construction}
\begin{lemma}\label{lem:Dd_normals}
The proper normal subgroups of $D_d = \langle r, s \mid r^d=s^2=1, srs=r^{-1} \rangle$ are:
\begin{enumerate}
    \item If $d$ is odd: precisely the subgroups of $\langle r \rangle$.
    \item If $d$ is even: the subgroups of $\langle r \rangle$, together with $\langle r^2, s \rangle$ and $\langle r^2, rs \rangle$.
\end{enumerate}
In the even case, $N_1 = \langle r^2, rs \rangle$ is transitive of index 2, while $N_2 = \langle r^2, s \rangle$ is intransitive of index 2.
\end{lemma}

\begin{lemma}\label{lem:decomp_lift}
Let $f = f_1 \circ \cdots \circ f_n$ be a composition of indecomposable rational functions. Let $H_n = \Gal(M_n/M_{n-1})$, where $M_i$ is the Galois closure of $f_1 \circ \cdots \circ f_i$ over $k(t)$. Then $H_n$ has the average fixed-point lifting property unless $M_n = M_{n-1}$.
\end{lemma}
\begin{proof}
The kernel of the restriction $H_n \to H_{n-1}$ factors through a product of primitive groups acting on subtrees. As every non-trivial normal subgroup of a primitive group is transitive, each factor acts transitively or trivially. Transitivity implies the lifting property via Lemma~\ref{lem:coset_burnside}; triviality implies $M_n = M_{n-1}$.
\end{proof}

\subsection{Proof of Theorem~\ref{thm:main_nonmart}}
The key geometric insight is that certain compositions yield unramified extensions of genus-zero function fields, which must be trivial.

\begin{lemma}[Unramified Cover Criterion]\label{lem:unram}
Let $f, g \colon \bP^1 \to \bP^1$ be rational functions. Let $K_1$ be the Galois closure of $f(x)-t$, and $L_2$ the Galois closure of $(f \circ g)(x)-t$. Suppose:
\begin{enumerate}
    \item $f(\CV(g)) \subseteq \CV(f)$,
    \item $\CV(g) \cap \Crit(f) = \emptyset$,
    \item For every $c \in \CV(f)$, if $e_c = \lcm\{e_z : f(z)=c\}$, then for all $y \in (f \circ g)^{-1}(c)$, $e_y$ divides $e_c$.
\end{enumerate}
Then $L_2/K_1$ is unramified.
\end{lemma}
\begin{proof}
Condition (1) ensures $L_2$ and $K_1$ share the same ramified primes. Conditions (2) and (3), combined with Abhyankar's Lemma (Lemma~\ref{lem:abhyankar}), imply ramification indices coincide. Hence $L_2/K_1$ is unramified.
\end{proof}

\begin{proof}[Proof of Theorem~\ref{thm:main_nonmart}]
Let $d \geq 4$ be even. Choose permutations $\sigma_1 = (1\,2\,\dots\,d)$, $\sigma_2 = (1\,2)(3\,4)\cdots$ a type 1 symmetry (i.e., with no fixed point), $\sigma_3 = \sigma_2 \sigma_1^{-1}$ a type 2 symmetry. Since $\sigma_1 \sigma_2 \sigma_3 = \operatorname{id}$ and $\langle \sigma_1, \sigma_2 \rangle \cong D_d$, the Riemann Existence Theorem \cite{Volklein} guarantees a cover $\varphi$ with ramification types $[d]$, $[2,\dots,2]$, $[2,\dots,2,1,1]$. By Riemann-Hurwitz, it has genus 0, hence a map $\phi: \bP^1 \rightarrow \bP^1$. Let $M/F$ be the field extension associated with $\phi$. By Lemma \ref{lem:abhyankar}, we see that the Galois closure of $M/F$ has ramification type $[m, m], [2,2,\dots,2], [2,2,\dots,2]$, which is the same as the one of $L/F$. So far we have showed that $\phi$ corresponds to an imprimitive Galois extension $L/F$ with a proper subfield $M/F$ of index 2. By counting degrees, we see that $\phi$ can be decomposed as $g \circ h$, where $g$ is of degree 2 which ramifies in the branches of $[m]$ and $[2,2,\dots,2]$ while splits in the branch of $[2,\dots,2,1,1]$. Since there always exists a linear fractional that sends the critical values of the branches $[m]$ and $[2,2,\dots,2]$ to the two unramified points in the branch of $[2,\dots,2,1,1]$, we may assume that $\phi = g \circ h$ has this property. In fact, we get a 1-parameter family of rational functions $\phi$ with this property as we have one variable of freedom from the linear fractional.

Let $K_1$ be the Galois closure of $\varphi(x)-t$, and $L_3$ the Galois closure of $(\varphi \circ g)(x)-t$. By Lemma~\ref{lem:unram}, $L_3/K_1$ is unramified. As a genus 0 covering has no non-trivial unramified extension, $L_3 = K_1$. The restriction map $\Gal(L_3/k(t)) \to \Gal(K_1/k(t))$ has trivial kernel, so the kernel $H_2 \leq G_2$ acts trivially on level-2 subtrees. This violates the lifting property \eqref{eq:lifting}. 
\end{proof}

\subsection{A Dihedral Example}
We first illustrate the construction in the simplest non-trivial case, $d=4$. Consider the polynomial $f(x) = x^4 - x^2 + 1$. It admits a decomposition into two quadratic polynomials:
\[ f = g_1 \circ h, \quad \text{where } g_1(x) = x^2 + x + 1 \text{ and } h(x) = x^2 - 1. \]
The critical points of $f$ are $\{0, \pm \frac{\sqrt{2}}{2}, \infty\}$. Analyzing the ramification behavior over their images:
\begin{itemize}
    \item $\infty$ is totally ramified with index $4$.
    \item The points $\pm \frac{\sqrt{2}}{2}$ map to $\frac{3}{4}$. Each has ramification type $(2,2)$ (two preimages of index $2$).
    \item The point $0$ maps to $1$. It has ramification type $(2,1,1)$ (one preimage of index $2$, two unramified preimages $\pm 1$).
\end{itemize}
To align the critical values for composition, we post-compose $f$ with a Möbius transformation $r \in \operatorname{PGL}_2(k)$ chosen so that
\[ r\left(f\left(\pm \tfrac{\sqrt{2}}{2}\right)\right) = -1, \qquad r(\infty) = 1. \]
Explicitly, for parameters $a \neq 0$ and $3a \neq -4b$, let
\[ r(x) = \frac{ax + b}{ax - \frac{3}{2}a - b}. \]
Setting $\varphi = r \circ f$ yields a one-parameter family of degree 4 rational functions:
\[ \varphi(x) = \frac{ax^4 - ax^2 + a + b}{ax^4 - ax^2 - a^2 - b}. \]
This family is invariant under simultaneous scaling of $(a,b)$, so it genuinely has one free parameter. Writing $g = r \circ g_1$, we verify the decomposition $\varphi = g \circ h$ where
\[ g(x) = \frac{ax^2 + ax + a + b}{ax^2 + ax - a^2 - b}. \]
By specializing parameters, we obtain infinitely many post-critically finite (PCF) maps. For instance, imposing $a + b = 0$ fixes the critical point $0$ in the forward orbit, yielding
\[ \psi(x) = \frac{2x^4 - 2x^2}{2x^4 - 2x^2 + 1}, \]
which has finite post-critical orbit $\operatorname{Orb}_\psi(0) = \{-1, 0\}$.

We now analyze the ramification in the Galois closure. Let $K_1$ be the splitting field of $\varphi(x) - t$ over $\overline{k}(t)$, and $L_2$ its Galois closure. The geometric monodromy group is $D_4$. Let $\pm u_1, \pm u_2$ denote the roots of $\varphi(x) - t$ in $L_2$. The intermediate fields $\overline{k}(u_1)/\overline{k}(t)$ and $\overline{k}(u_2)/\overline{k}(t)$ both have degree $4$. 

At the prime corresponding to $\infty$, the place is totally ramified (index 4) in both intermediate extensions. Since $L_2$ is the compositum of conjugate fields, Abhyankar's Lemma (Lemma~\ref{lem:abhyankar}) implies the ramification index in $L_2$ is $\lcm(4,4) = 4$ at each of the two primes above $\infty$. Thus $\infty$ has ramification type $[4,4]$ in $L_2$. Similarly, the critical values at $\pm \frac{\sqrt{2}}{2}$ have ramification index 2 in each intermediate field, yielding type $[2,2,2,2]$ in $L_2$. By degree counting and Galois transitivity, the remaining critical value (from $0$) must also have type $[2,2,2,2]$ in $L_2$.

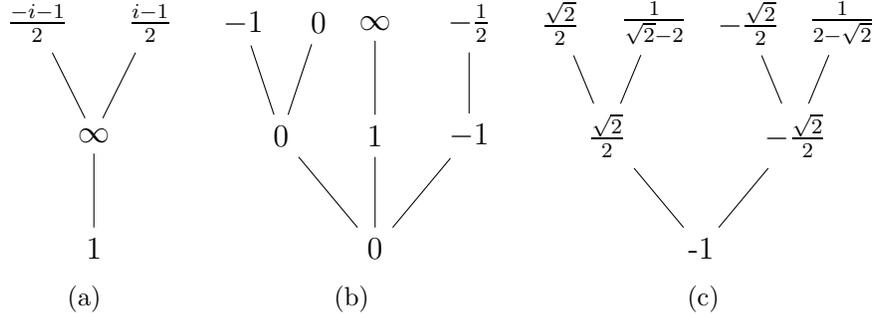
\begin{figure}
\centering
\subfigure[]{
\begin{tikzpicture}[level distance= 1.5cm,
level 1/.style={sibling distance= 2cm},
level 2/.style={sibling distance= 1.5cm}, level 3/.style={sibling distance= 1.5cm}]

    \node (Root) {1} [grow=up]
   child { node {$\infty$}
            child{ node {$\frac{i-1}{2}$}}
            child{ node {$\frac{-i-1}{2}$}}
        }
;
\end{tikzpicture}
}
\subfigure[]{
\begin{tikzpicture}[level distance= 1.5cm,
level 1/.style={sibling distance= 1.25cm},
level 2/.style={sibling distance= 1cm}, level 3/.style={sibling distance= 1.5cm}, level 4/.style={sibling distance= 0.7cm}]
    \node (Root) {0} [grow=up]
    child { node {$-1$}
            child{ node {$-\frac{1}{2}$}}
    }
   child { node {$1$}
            child{ node {$\infty$}}
        }
   child { node {$0$}
            child{ node {$0$}}
            child{ node {$-1$}}
   }
;
\end{tikzpicture}
}
\subfigure[]{
\begin{tikzpicture}[level distance= 1.5cm,
level 1/.style={sibling distance= 2.5cm},
level 2/.style={sibling distance= 1.25cm}, level 3/.style={sibling distance= 1.5cm}, level 4/.style={sibling distance= 0.7cm}]

    \node (Root) {-1} [grow=up]
        child { node {$-\frac{\sqrt{2}}{2}$}
            child{ node {$\frac{1}{2 - \sqrt{2}}$}}
            child{ node {$-\frac{\sqrt{2}}{2}$}}
        }
	    child { node {$\frac{\sqrt{2}}{2}$}
            child{ node {$\frac{1}{\sqrt{2} - 2}$}}
            child{ node {$\frac{\sqrt{2}}{2}$}}
    }
;
\end{tikzpicture}
}
\caption{Ramification Portraits of $g \circ \psi$.}
\end{figure}
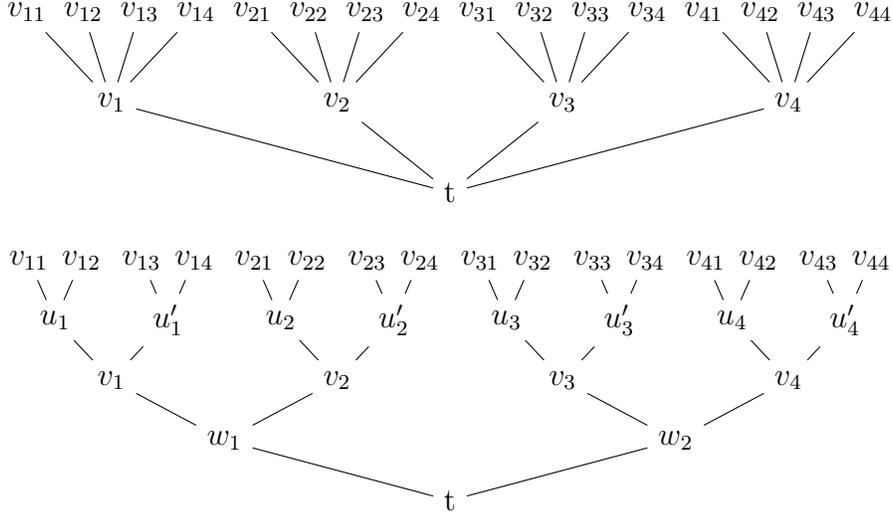
\begin{figure}
\centering
\begin{tikzpicture}[level distance= 0.8cm,
level 1/.style={sibling distance= 6cm},
level 2/.style={sibling distance= 3cm}, level 3/.style={sibling distance= 1.5cm}, level 4/.style={sibling distance= 0.7cm}]

    \node (Root) {t} [grow=up]
    child {
    node {$w_2$} 
    child { node {$v_4$} 
        child {node {$u_4'$}
            child {node {$v_{44}$}}
            child {node {$v_{43}$}}
            }
        child {node {$u_4$}
            child {node {$v_{42}$}}
            child {node {$v_{41}$}}
        }
        }
    child { node {$v_3$} 
        child {node {$u_3'$}
            child {node {$v_{34}$}}
            child {node {$v_{33}$}}
        }
        child {node {$u_3$}
            child {node {$v_{32}$}}
            child {node {$v_{31}$}}
        }
        }
}
child {
    node {$w_1$}
    child { node {$v_2$} 
        child {node {$u_2'$}
            child {node {$v_{24}$}}
            child {node {$v_{23}$}}
        }
        child {node {$u_2$}
            child {node {$v_{22}$}}
            child {node {$v_{21}$}}
        }
        }
    child { node {$v_1$} 
        child {node{$u_1'$}
            child {node {$v_{14}$}}
            child {node {$v_{13}$}}
        }
        child {node {$u_1$}
            child {node {$v_{12}$}}
            child {node {$v_{11}$}}
        }
    }
};
\end{tikzpicture}
\caption{We get a level $4$ tree by decomposing $\phi$ as $g \circ h$.}
\end{figure}

Applying the Riemann--Hurwitz formula to the cover $C \to \bP^1_{\overline{k}}$ corresponding to $K_1/\overline{k}(t)$:
\[ g(C) = 1 + \frac{1}{2} \left( -2 \cdot 4 + \sum_{P \in C} (e_P - 1) \right) = 1 + \frac{1}{2}(-8 + 4+2+2+2+2) = 0. \]
Thus $K_1 = \overline{k}(C)$ is a genus 0 function field.

The composition structure $\varphi = g \circ h$ induces an isomorphism between the 4-ary tree of level 2 for $\varphi$ and the binary tree of level 4 for $h$. Computing the ramification portrait of $\varphi \circ g$ explicitly, we find that its Galois closure $L_3$ has exactly the same ramification type as $K_1$: one branch of type $[4,4]$ and two branches of type $[2,2,2,2]$. Since both $L_3$ and $K_1$ are genus 0 function fields and $L_3/K_1$ is an extension, the Riemann--Hurwitz formula forces $g(L_3) = g(K_1) = 0$. Since a genus 0 cover admits no non-trivial unramified extensions, $L_3 = K_1$.

Consequently, the restriction map $\pi_2 \colon \Gal(L_3/\overline{k}(t)) \to \Gal(K_1/\overline{k}(t))$ is an isomorphism, and its kernel $H_2 = \ker(\pi_2)$ is trivial. This means $H_2$ acts trivially on each level-2 subtree above a fixed vertex in $L_1$. In particular, it acts intransitively on the subtrees $u_1^*$, violating the average fixed-point lifting property. Hence the fixed-point process for $\varphi$ is not a martingale.

\subsection{Characterization in Genus Zero}
Having established existence, we ask whether non-martingale behavior can arise from other monodromy groups. The essential observation is that the martingale failure requires $G_1$ to contain a normal transitive subgroup and a normal intransitive subgroup of identical index. Under the genus-zero assumption, the following Klein's classification restricts $G_1$ to a short list. 
\begin{lemma}[Klein's Classification]\label{lem:klein}
The finite subgroups of $\operatorname{PGL}_2(\overline{k})$ are isomorphic to one of the following:
\begin{enumerate}
    \item The cyclic group $C_n$ for $n \geq 1$.
    \item The dihedral group $D_n$ of order $2n$ for $n \geq 2$.
    \item The tetrahedral group $A_4$.
    \item The octahedral group $S_4$.
    \item The icosahedral group $A_5$.
\end{enumerate}
\end{lemma}
\begin{proof}
See \cite{Volklein} or standard references on the geometry of $\bP^1$.
\end{proof}
We now prove that even dihedral groups are the only obstruction to the martingale property under the genus-zero assumption.
\begin{theorem}\label{thm:char_even_dihedral}
Assume that the geometric Galois closure $L_1$ of $K_1/\overline{k}(t)$ has genus $0$. If the fixed-point process of $f$ is not a martingale, then the geometric monodromy group $G_1 = \Gal(K_1/\overline{k}(t)) \cong D_m$ for some even integer $m$.
\end{theorem}

\begin{proof}
Assume the fixed-point process is not a martingale. By Lemma~\ref{lem:mart_equiv} and Theorem~\ref{thm:bjkl_mart}, this implies that the group $G_1$ does not satisfy the average fixed-point lifting property at the second level. Consequently, $G_1$ must possess a normal transitive subgroup $N_1$ and a normal intransitive subgroup $N_2$ of the same index $p$.

Since $L_1$ has genus $0$, the Galois closure of the extension corresponds to a curve $C$ of genus $0$, and the Galois group $G_1$ embeds into $\Aut(C) \cong \operatorname{PGL}_2(\overline{k})$. By Lemma~\ref{lem:klein}, $G_1$ must be one of $C_n$, $D_n$, $A_4$, $S_4$, or $A_5$. We exclude the possibilities one by one.

\textbf{Case 1: $G_1 \cong D_n$ with $n$ odd.}
We show that the process is a martingale when $n$ is odd. The group $D_n = \langle r, s \mid r^n=s^2=1, srs=r^{-1} \rangle$ acts on the $n$ roots.
For any $i \in \{1, \dots, n\}$, the stabilizer $H_i = \Stab_{D_n}(i)$ consists of the identity $e$ and exactly one reflection $s_i$. Thus $|H_i| = 2$.
Let $K$ be the kernel of the restriction map to the first level (the subgroup acting trivially on $X_1$). Since $G$ is recurrent, the section map $\phi_i : \Stab_G(i) \to \Aut(i^*) \cong G$ is surjective onto $G_1 \cong D_n$.
The restriction of $K$ to the subtree $i^*$, denoted $\phi_i(K)$, is a normal subgroup of $D_n$. Since $|H_i| = 2$, the index $[\phi_i(H_i) : \phi_i(K)]$ must divide 2.
The only normal subgroups of $D_n$ with index dividing 2 are $D_n$ itself (index 1) and the rotation subgroup $\langle r \rangle$ (index 2). Since $n$ is odd, both $D_n$ and $\langle r \rangle$ act transitively on the $n$ points of the subtree $i^*$. Therefore, the stabilizer section acts transitively on the subtree for all $i$, satisfying the martingale criterion. Thus, $D_n$ for $n$ odd always yields a martingale process.

\textbf{Case 2: $G_1 \in \{A_4, S_4, A_5\}$.}
We check if these groups possess a pair of normal subgroups $(N_1, N_2)$ with $N_1$ transitive, $N_2$ intransitive, and $|G/N_1| = |G/N_2|$.
\begin{itemize}
    \item \emph{$A_5$:} $A_5$ is simple, so it has no proper non-trivial normal subgroups.
    \item \emph{$A_4$:} The only proper non-trivial normal subgroup is the Klein four-group $V_4$. $V_4$ acts transitively on 4 points (as it contains all double transpositions). There are no other normal subgroups, so no intransitive candidate $N_2$ exists.
    \item \emph{$S_4$:} The normal subgroups are $A_4$ (index 2, transitive) and $V_4$ (index 6, intransitive). The indices $2$ and $6$ do not match. Thus no suitable pair exists.
\end{itemize}

\textbf{Case 3: $G_1 \cong C_n$.}
A cyclic group $C_n$ acting transitively on $n$ points is generated by an $n$-cycle. Any subgroup $H \leq C_n$ is normal. $H$ acts transitively if and only if $H = C_n$. If $H$ is a proper subgroup, it is intransitive. However, we require a proper transitive subgroup $N_1$ to satisfy the condition of Lemma~\ref{lem:mart_equiv} (specifically, $N_1$ must be transitive and proper to allow for a non-trivial intransitive quotient comparison). Since $C_n$ is the smallest transitive group of degree $n$, it has no proper transitive subgroup. Thus, $C_n$ cannot support a non-martingale structure.

\textbf{Case 4: $G_1 \cong D_n$ with $n$ even.}
By Lemma~\ref{lem:Dd_normals}, $D_n$ possesses the normal subgroup $N_1 = \langle r^2, rs \rangle$ which is transitive of index 2, and the normal subgroup $N_2 = \langle r^2, s \rangle$ which is intransitive of index 2. This pair satisfies the necessary and sufficient conditions for the failure of the martingale property.

Therefore, the only remaining possibility is $G_1 \cong D_n$ with $n$ even.
\end{proof}

\begin{remark}
Our construction explicitly provides rational functions with $g(L_1) = 0$ demonstrating this phenomenon. If the genus-zero assumption is dropped, there may exist rational functions with non-martingale fixed-point processes associated with monodromy groups of higher-genus covers. Proving such existence would require new techniques, which lies beyond the scope of this paper.
\end{remark}

We conclude this section with a question linking the arithmetic structure of the rational function to the group-theoretic pattern of its iterated monodromy group.

\begin{question}\label{q:finite_type}
Let $f\in k(x)$ be a non-PCF rational function of degree $d\ge 2$.
Is the profinite iterated monodromy group $G(f)$ necessarily a recurrent group of finite type generated by a pattern group $P$ of depth $m$, where $m$ depends on $f$?
\end{question}

This is known to hold for PCF maps~\cite{similar, pink2013infinite} and for generic polynomials~\cite{Odoni}.
However, for non-PCF maps the iterated monodromy group can be infinitely generated~\cite{Bar16}, so the answer may be negative.
A positive answer would imply that non-martingale behaviour is controlled by finite ramification patterns; a negative answer would show that deeper, non-finite-type structures can also produce martingale failure.
We leave this as an open problem for future investigation.

\section{Martingale Fixed-Point Processes in Arbitrary Dimensions}\label{sec:higher_dim}

A natural question is whether the martingale criterion extends to higher-dimensional morphisms, where the wreath product structure $G_n \cong [G_1]^n$ is known to fail in general (cf. Example~\ref{ex:split}). 

We show that the wreath product hypothesis can be replaced by an orbit-disjointness condition. The key insight is that transitivity on level-one subtrees depends only on the inertia groups lifting unramified primes in the relative extension $K_n/K_{n-1}$. When critical orbits remain disjoint, these inertia groups lift faithfully and generate transitive actions on each fiber, preserving the martingale property without requiring a wreath product decomposition. 

\subsection{Geometric Setup and the Split Morphism Counterexample}
Let $X \subset \bP^N_k$ be a normal variety of dimension $N \geq 1$ over an algebraically closed field $k$. Let $f \colon X \to X$ be a finite, generically \'etale morphism of degree $d$. The induced pullback $f^* \colon K(X) \to K(X)$ is a degree $d$ extension of function fields. For each $n \geq 1$, define $K_n$ to be the Galois closure of the extension induced by $(f^n)^*$ over $K(X)$, and let $G_n = \Gal(K_n/K(X))$. Set $H_n = \ker(G_n \to G_{n-1})$ for $n \geq 2$, with $H_1 = G_1$.

In contrast to dimension one, the iterated Galois groups $G_n$ need not be wreath powers. The following example demonstrates this failure explicitly.

\begin{example}\label{ex:split}
Let $\Phi \colon \bP^1 \times \bP^1 \to \bP^1 \times \bP^1$ be a split morphism given by $\Phi(x,y) = (\Phi_1(x), \Phi_2(y))$, where $\Phi_1, \Phi_2$ are disjoint rational functions. The splitting field $M$ of $\Phi(x,y) = (t_1, t_2)$ over $k(t_1, t_2)$ is the compositum of the splitting fields $L_1$ of $\Phi_1(x)=t_1$ and $L_2$ of $\Phi_2(y)=t_2$. Since the fields are linearly disjoint, we have
\[ \Gal(M/k(t_1, t_2)) \cong \Gal(L_1/k(t_1)) \times \Gal(L_2/k(t_2)) =: H_1 \times H_2. \]
Iterating yields $G_n \cong [H_1]^n \times [H_2]^n$. Unless one of $H_1, H_2$ is trivial, $G_n$ is not isomorphic to the wreath power $[H_1 \times H_2]^n$. Thus, the strict wreath product condition of \cite[Theorem 3.1]{wreath} fails for generic split maps, even when the fixed-point process remains a martingale.
\end{example}

To handle this, we replace the standard preimage tree with a tree of cosets. Let $A_0 = \{G_1\}$ and for $n \geq 1$, let $A_n$ be the set of left cosets of $H_{n+1}$ in $G_{n+1}$. Then $\#A_n = d^n$, and $G_{n+1}$ acts transitively on $A_n$ by left multiplication. This coset tree is isomorphic to the regular $d$-ary tree $X_n$, and $G_n$ acts on it via restriction of the $G_{n+1}$-action. The fixed-point process $\{Y_n\}$ is defined identically: $Y_n(g) = \#\Fix(\pi_n(g))$ with respect to this action.

\subsection{The Higher-Dimensional Martingale Criterion}
Let $\Delta_f \subset X$ denote the critical divisor of $f$, i.e., the set of prime divisors in $X$ where the ramification index exceeds 1. We establish conditions under which the fixed-point process remains a martingale despite potential failure of the wreath product structure.

\begin{lemma}\label{lem:ramified_primes}
Let $n < N$. Every prime of $X$ that is ramified under $f^n$ is of the form $f^j(q)$ for some prime $q \in \Delta_f$ and some integer $j \leq n$.
\end{lemma}
\begin{proof}
The base case $n=1$ holds by definition. Assume the claim holds for $n-1$. Consider the field tower $K(X) \subset K_{n-1} \subset K_n$. Multiplicativity of ramification indices gives $e(\mathfrak{P}|\mathfrak{p}) = e(\mathfrak{P}|\mathfrak{q}) e(\mathfrak{q}|\mathfrak{p})$, where $\mathfrak{q} = \mathfrak{P} \cap K_{n-1}$. A prime $\mathfrak{P}$ ramifies in $K_n/K(X)$ only if it ramifies in $K_{n-1}/K(X)$ or if $\mathfrak{q}$ ramifies in $K_n/K_{n-1}$. By the inductive hypothesis, ramification in $K_{n-1}$ occurs only at primes of the form $f^j(q)$ with $j \leq n-1$. Ramification in the relative extension $K_n/K_{n-1}$ occurs precisely above primes mapping to $\Delta_f$ under $f$. Hence any ramified prime in $K_n$ must be of the form $f(f^j(q)) = f^{j+1}(q)$ with $j+1 \leq n$, completing the induction.
\end{proof}

\begin{theorem}\label{thm:martingale_high_dim}
Suppose there exists a subset $S \subseteq \Delta_f$ such that:
\begin{enumerate}
    \item For all $p \in S$ and $q \in \Delta_f$, and all integers $m, r \geq 0$, we have $f^m(p) \neq f^r(q)$ unless either $(p=q \text{ and } m=r)$ or $(p \neq q \text{ and } m \leq r)$.
    \item The group $G_1$ is generated by the inertia groups 
    $$ \{ I(\mathfrak{p}|f(p)) : p \in S, \mathfrak{p} \text{ lies over } f(p) \text{ in } K_1 \}. $$
\end{enumerate}
Then the fixed-point process associated to $G(f)$ is a martingale.
\end{theorem}

\begin{proof}
By Theorem~\ref{thm:bjkl_mart}, it suffices to prove that for every $n \geq 1$, the kernel $H_{n+1} = \ker(G_{n+1} \to G_n)$ acts transitively on each level-one subtree $\omega^*$ rooted at $\omega \in A_n$. 

Fix $\omega \in A_n$. The stabilizer $\Stab_{G_{n+1}}(\omega)$ projects isomorphically onto $G_1$ via the restriction map since elements fixing $\omega$ act faithfully on the $d$
preimages of $\omega$ under $f$, yielding a subgroup $G' \leq G_{n+1}$ with $G' \cong G_1$ that acts on the subtree $\omega^*$ of size $d$. Let $L'_\omega / K_n$ be the field extension corresponding to the Galois closure of this subtree action. By Condition (2), $G'$ is generated by inertia groups over primes lying above the critical values $f(p)$ for $p \in S$.

To apply Abhyankar's Lemma to the compositum tower $K_{n+1} = K_n \cdot L'_\omega$ over $K_n$, we must ensure the branch locus of the relative extension is unramified in the base tower. In the relative extension $L'_\omega/K_n$, the branch points correspond to $f(p)$ for $p \in S$. The projection of these points down to the base variety $X$ via $f^n$ yields exactly $f^{n+1}(p)$. By Lemma~\ref{lem:ramified_primes}, the branch locus of $K_n/K(X)$ consists of forward images of critical values up to level $n-1$, which are points of the form $f^r(f(q)) = f^{r+1}(q)$ for $q \in \Delta_f$ and $0 \leq r \leq n-1$. 

We must show that $f^{n+1}(p) \neq f^{r+1}(q)$. Setting $m = n+1$ and $k = r+1$, we are comparing $f^m(p)$ and $f^k(q)$ where $1 \leq k \leq n$. Since $m = n+1 > n \geq k$, we strictly have $m > k$. Condition (1) explicitly prohibits $f^m(p) = f^k(q)$ when $m > k$. Geometrically, this means the forward critical orbits originating from $S$ do not intersect the branch locus of $K_n/K(X)$. Consequently, the primes in $K_n$ ramifying in $L'_\omega$ are unramified in $K_n/K(X)$.

Let $\mathfrak{P}''$ be a prime in $K_{n+1}$ extending such a branch point $f^{n+1}(p)$, let $\mathfrak{q}_n$ be the prime below it in $K_n$, and let $\mathfrak{P}_1$ be the prime below it in $L'_\omega$. Because the underlying prime $f^{n+1}(p)$ is unramified in $K_n$, its ramification index is $e(\mathfrak{q}_n | f^{n+1}(p)) = 1$. Thus, the restriction map
\[ I(\mathfrak{P}''|f^{n+1}(p)) \longrightarrow I(\mathfrak{P}_1|f(p)) \]
is surjective. 

Now apply Lemma~\ref{lem:inertia}: since the prime is unramified in the intermediate field $K_n$, the inertia group $I(\mathfrak{P}''|f^{n+1}(p))$ acts trivially on $K_n$, yielding an isomorphism
\[ I(\mathfrak{P}''|f^{n+1}(p)) \cong I(\mathfrak{P}''|\mathfrak{q}_n). \]
The group $I(\mathfrak{P}''|\mathfrak{q}_n)$ is contained in $H_{n+1}$ and generates the inertia subgroup of $H_{n+1}$ corresponding to the branch point $f^{n+1}(p)$. Since this holds for every branch point arising from $S$, and $G_1$ is generated by these inertia groups by Condition (2), it follows that the relative group $G' \cong G_1$ surjects onto the action of $H_{n+1}$ on $\omega^*$. 

Because $G_1$ acts transitively on $A_1$, every restriction $G'|_{\omega^*}$ acts transitively on the $d$ vertices of $\omega^*$. Therefore, $H_{n+1}$ acts transitively on each level-one subtree. By Theorem~\ref{thm:bjkl_mart}, the fixed-point process is a martingale.
\end{proof}

\begin{corollary}
Under the assumptions of Theorem~\ref{thm:martingale_high_dim}, the sequence $\{Y_n\}$ converges almost surely with respect to the Haar measure on $G(f)$.
\end{corollary}

\begin{remark}
Theorem~\ref{thm:martingale_high_dim} is inspired by the wreath product criterion of \cite[Theorem 3.1]{wreath}. Condition (1) permits controlled overlap in critical orbits (e.g., $p \neq q$ with $m \leq r$), which occurs naturally in split morphisms and certain PCF maps, while still guaranteeing transitivity on subtrees. For example, Chebyshev polynomials in higher dimensions satisfy these inertia conditions but yield $G_n \not\cong [G_1]^n$ \cite[Proposition 1.2]{martingale}; On the other hand, $H_{n}$ is generated by the ramification group over the point at infinity, which has degree $d^n$ and acts transitively on each subtree at the $n$-th level. That means it can induce a martingale as above.
\end{remark}

\section*{Acknowledgements}
This article is based on a chapter of the second author's doctoral thesis
completed at the University of Rochester, 2025. The authors are extremely grateful to their PhD. advisor Thomas Tucker for guiding us into this fascinating arithmetic world. We are deeply indebted to his enlightening discussions during the preparation of this paper and proofreading the first draft. The second author would also like to thank Juan Rivera-Letelier and Dinesh Thakur for their wonderful lectures that triggered his interest in arithmetic geometry and dynamical systems. 

\bibliographystyle{amsalpha}
\bibliography{Bib}
\end{document}